\documentclass[12pt,twoside]{amsart}

\usepackage{amsfonts}
\usepackage{amsmath}
\usepackage{epsfig}
\usepackage{latexsym}

\newtheorem{thm}{{\sc Theorem}}

\newtheorem{cor}{{\sc Corollary}}
\newtheorem{defn}{{\sc Definition}}

\title{Equidimensionality of  complex Lagrangian fibrations}
\author{Daisuke Matsushita}
\subjclass{Primary 14E35, Secondary 14D05}
\address{Research Institue for Mathematical Sciences \\
         Kyoto University, Oiwake-Cyo Kitashirakawa \\
         Sakyo-Ku Kyoto 606-8052 Japan}
\thanks{*Research Fellow of the Japan Society for the Promotion of Science} 
\email{tyler@kurims.kyoto-u.ac.jp}

\begin{document}
\maketitle
\begin{abstract}
 We prove that every irreducible component of a fibre 
 of a complex Lagrangian fibration is Lagrangian subvariety.
 Especially, complex Lagrangian fibations are equidimensional. 
\end{abstract}
\section{Introduction}
 First we define {\it Lagrangian subvariety \/}.
\begin{defn}
 Let $X$ be a manifold with a holomorhpic symplectic form $\omega$.
 A subvariety $Y$ is said to be  Lagrangian subvariety
 if $\dim Y = (1/2)\dim X$ and
 there exists a resolution $\nu : \tilde{Y} \to Y$ such that
 $\nu^{*}\omega $ is identically zero on $\tilde{Y}$.
\end{defn}

\noindent
 Note that this notion does not depend the choice of $\nu$.
 We prove the following theorem.
\begin{thm}\label{main}
  Let $f: X \to B$ be a proper surjective morphism 
  over a normal base $B$.
  Assume that $X$ is a K\"{a}hler manifold with
  holomorphic symplectic form $\omega$ and a general fiber
  of $f$ is a Lagrangian submanifold with respect to $\omega$.
  Then every irreducible component of a fibre of $f$ is
  a Lagrangian subvariety. 
\end{thm}
\noindent
  From Theorem \ref{main}, we obtain the following corollary.
\begin{cor}\label{equi_dimensional}
  Under the situation above, $f$ is equidimensional. 
  Especially $f$ is flat if $B$ is smooth.
\end{cor}

\noindent 
 Combining Corollary \ref{equi_dimensional} 
 with \cite[Theorem 2]{matsu} and
 \cite[Theorem 1]{matsu2}, we obtain the following result.
\begin{cor}
  Let $f:X \to B$ be a fibre space of an irreducible
  symplectic manifold $X$ over a normal
  K\"{a}hler base $B$. Then $f$ is equidimensional.
\end{cor}

\vspace{5mm}
\noindent
{\sc Acknowledgment. \quad} 
 The author express his thanks to Professors 
 A.~Beauville and A.~Fujiki for their advice and encouragement. 

\section{Proof of Theorem \ref{main}}
  We refer the following theorem due to Koll\'{a}r 
 \cite[Theorem 2.2]{kollar} and Saito \cite[Theorem 2.3, Remark 2.9]{saito} 
\begin{thm}\label{kollar}
  Let $f: X \to B$ be a proper surjective morhpism from
  a smooth K\"{a}hler manifold $X$ to a normal variety $B$.
  Then ${\rm R}^{i}f_{*}\omega_{X}$ is torsion free, where $\omega_X$ 
  is the dualizing sheaf of $X$.
\end{thm}
\noindent
  Let $\bar{\omega}$ be the complex conjugate of $\omega$.
  By Leray spectral sequence, there exists a morphism
$$
 H^{2}(X , {\mathcal O}_{X}) \to H^{0}(B, {\rm R}^{2}f_{*}{\mathcal O}_{X}).
$$
  Then $\bar{\omega}$ is a torsion element in 
  $H^{0}(B, {\rm R}^{2}f_{*}{\mathcal O}_{X})$ since
  a general fibre of $f$ is a Lagrangian manifold.
  In addition, $\omega_{X} \cong {\mathcal O}_{X}$.
  Hence
  $\bar{\omega}$ is zero in 
  $H^{0}(B, {\rm R}^{2}f_{*}{\mathcal O}_{X})$ by Theorem \ref{kollar}.
  We derive a contradiction assuming 
  that there exists an irrducible component of a fibre of $f$
  which is not a Lagrangian subvariety.
  The letter $V$ denotes an non Lagrangian component.
  We take an embedding resolution $\pi : \tilde{X} \to X$ of 
  $V$. Let $\tilde{V}$ be the proper transform of $V$.
  We will show that $\pi^{*}\omega$ is not zero in
  $H^{0}(\tilde{V},\Omega^{2}_{\tilde{V}})$.
  If $\dim V = (1/2)\dim X$, it is obious by the definition.
  If $\dim V > (1/2)\dim X$, we take a smooth point $q \in V$
  such that $\pi$ is isomorhpic in a neighborhood of $q$.
  Since $\dim V > (1/2)\dim X$ and $\omega$ is nondegenerate,
  the restriction of $\omega$ on 
  the tangent space of $V$ at $q$ is nonzero.
  Because $\pi$ is isomorphic in a neighborhood of $q$,
  $\pi^{*}\omega$ is not zero in
  $H^{0}(\tilde{V},\Omega^{2}_{\tilde{V}})$.
  Take the complex conjugate, $\pi^{*}\bar{\omega}$
  is not zero in $H^{2}(\tilde{V}, {\mathcal O}_{\tilde{V}})$.
  Therefore $\bar{\omega}$ is not zero in 
  $H^{2}(V,{\mathcal O}_V)$. 
  Let $p := f(V)$ and $X_p := f^{-1}(p)$.
  We consider the following morphism:
$$
  {\rm R}^{2}f_{*}{\mathcal O}_{X}\otimes k(p)
  \to
  H^{2}(X_p ,{\mathcal O}_{X_p})
  \to
  H^{2}(V, {\mathcal O}_{V}).
$$
  Then $\bar{\omega}$ is zero in 
  ${\rm R}^{2}f_{*}{\mathcal O}_{X}\otimes k(p)$
  and notzero in $H^{2}(V, {\mathcal O}_{V})$.
  That is a contradiction.
\qed

\end{document}